\magnification=\magstep1
\overfullrule=0pt

\font\twelvebf=cmbx12 scaled 1000

\font\tenmsb=msbm10 scaled 1000
\font\sevenmsb=msbm7 scaled 1000
\font\fivemsb=msbm5 scaled 1000
\newfam\msbfam
\textfont\msbfam=\tenmsb
\scriptfont\msbfam=\sevenmsb
\scriptscriptfont\msbfam=\fivemsb
\def\Bbb#1{{\fam\msbfam\relax#1}}

\def\N{{\Bbb N}}
\def\E{{\Bbb E}}
\def\R{{\Bbb R}}
\def\Q{{\Bbb Q}}
\def\C{{\Bbb C}}
\def\cC{{\cal C}}
\def\cF{{\cal F}}
\def\cG{{\cal G}}
\def\cM{{\cal M}}
\def\cS{{\cal S}}
\def\cT{{\cal T}}
\def\cU{{\cal U}}
\def\cL{{\cal L}}

\def\ae{\quad\hbox{almost everywhere}}
\def\as{\quad\hbox{almost surely}}

\def\modo#1{{\left|#1\right|}}
\def\normo#1{{\left\|#1\right\|}}
\def\smodo#1{{\mathopen|#1\mathclose|}}

\def\QED{\par\hfill Q.E.D.}

\def\Proof:{\medskip\noindent{\bf Proof:}\ \ } 
\def\Proofof#1:{\medskip\noindent{\bf Proof of #1:}\ \ }

\def\moreproclaim{\par}

{
\nopagenumbers

\null\vfill
{
\centerline{\hbox{
\vrule height -0.4 pt depth 0.8 pt width 26.5 em
\kern - 26.5 em
\raise  0.03ex  \hbox{\bf E} 
\raise  0.06ex \hbox{l} 
\raise .13ex \hbox{e}
\raise .24ex \hbox{c} 
\raise .45ex \hbox{t} 
\raise .78ex \hbox{r} 
\raise 1.31ex \hbox{o} 
\raise 2.08ex \hbox{n} 
\raise 3.14ex \hbox{i} 
\raise 4.53ex \hbox{c} 
\kern   1em
\raise 8.15ex \hbox{\bf J} 
\raise 10.15ex \hbox{o} 
\raise 12.04ex \hbox{u} 
\raise 13.60ex \hbox{r} 
\raise 14.64ex \hbox{n} 
\kern .3 em
\vrule
\kern -.3em
\raise 15ex \hbox{a} 
\raise 14.64ex \hbox{l} 
\kern   1em
\raise 12.04ex \hbox{o} 
\raise 10.15ex \hbox{f} 
\kern   1em
\raise 6.23ex \hbox{\bf P} 
\raise 4.53ex \hbox{r} 
\raise 3.14ex \hbox{o} 
\raise 2.08ex \hbox{b} 
\raise 1.31ex \hbox{a} 
\raise .78ex \hbox{b} 
\raise .45ex \hbox{i} 
\raise .24ex \hbox{l} 
\raise .13ex \hbox{i} 
\raise .06ex \hbox{t} 
\raise .03ex \hbox{y}
}}}

\bigskip

\centerline{\tt http://www.math.washington.edu/\~{}ejpecp/}

\vskip 0.5 true in

\centerline{\twelvebf Concrete Representation of Martingales}

\bigskip

\centerline{\bf Stephen Montgomery-Smith}
\centerline{\it Department of Mathematics}
\centerline{\it University of Missouri, Columbia, MO 65211}
\centerline{\tt stephen@math.missouri.edu}
\centerline{\tt http://math.missouri.edu/\~{}stephen}

\bigskip
\centerline{Volume~3, paper number~15, 1998}
\centerline{\tt http://www.math.washington.edu/\~{}ejpecp/EjpVol3/paper15.abs.html}

\bigskip
\centerline{15 pages, submitted June 4, 1998, published December 2, 1998}

\bigskip
\noindent {\bf Abstract:}\ \ 
Let $(f_n)$ be a mean zero
vector valued martingale sequence.  Then there exist
vector valued functions $(d_n)$ from $[0,1]^n$ such that
$\int_0^1 d_n(x_1,\dots,x_n)\,dx_n = 0$ for almost all
$x_1,\dots,x_{n-1}$, and such that the law of $(f_n)$
is the same as the law of 
$(\sum_{k=1}^n d_k(x_1,\dots,x_k))$.
Similar results for tangent sequences and sequences satisfying condition~(C.I.)
are presented.  
We also present a weaker 
version of a result of McConnell that
provides a Skorohod like representation for vector valued
martingales.


\bigskip
\noindent A.M.S.\ Classification (1991): 60G42, 60H05.

\bigskip
\noindent Research supported in part by the N.S.F.\ and the Research Board of
the University of Missouri.

\vfill\eject
}
\pageno=2

\beginsection 1.\ \ Introduction

In this paper, we seek to give a concrete representation of martingales.
We will present theorems of the following form.  Given a martingale
sequence $(f_n)$ (possibly vector valued), 
there is a martingale $(g_n)$ on the probability
space $[0,1]^\N$, with respect to the filtration 
$\cL_n$, the minimal sigma field for which the first $n$ coordinates
of $[0,1]^\N$ are measurable, such that the sequence $(f_n)$
has the same law as $(g_n)$.  Thus any martingale may be 
represented
by a martingale
$$ g_n((x_n)) = \sum_{k=0}^n d_k(x_1,\dots,x_k) ,$$
where $\int_0^1 d_n(x_1,\dots,x_n) \, dx_n = 0$ for $n \ge 1$.  (Here, as
in the rest of this paper, the notation $(x_n)$ will refer to the
sequence $(x_n)_{n \in \N}$, where $\N$ refers to the positive integers.)

The value of such a result is perhaps purely psychological.  However,
we will also present similar such results for tangent sequences,
and also for sequences satisfying condition (C.I.).
Tangent sequences and condition (C.I.), as defined in this paper, 
were introduced by Kwapie\'n and Woyczy\'nski [KW].
The abstract definition of tangent sequences can be, perhaps, a
little hard to grasp.  However, in this concrete setting, it is
clear what their meaning is.  To demonstrate the psychological
advantage that this view gives, we will give a new proof
of a result of McConnell [M2] that states that tangent
sequence inequalities hold for the UMD spaces.

We will also give a weaker version of a result of McConnell [M1]
that provides a Skorohod representation theorem for
vector valued martingales, in that any martingale is a stopped
continuous time stochastic process.  
This means that many martingale inequalities that
are true for continuous time stochastic processes are automatically also 
true for
general martingales.

\beginsection 2.\ \ Representations of sequences of random variables

In this section, we present the basic result upon which all our other
results will depend.  

First let us motivate this result by considering just one random
variable, that is a measurable
function $f$ from the underlying probability space to a separable
measurable space $(S,\cS)$.
(A measurable space is said to be {\it separable\/} if its sigma field is
generated by a countable
collection of sets.)
We seek to find a measurable function
$g:[0,1] \to S$ that has the same law as $f$, that is, given any
measurable set $A$, we have that $\Pr(f \in A) = \lambda(g \in A)$,
or equivalently, given any measurable bounded function $F:S \to \R$, we have
that $E(F(f)) = \int F(g) \, d\lambda$.
Here $\lambda$ refers to the Lebesgue measure on $[0,1]$.

First, it may be shown without loss of generality that $(S,\cS)$
is $\R$ with the Borel sets.  (The argument for this may be found in
Chapter~1 of [DM], and is presented below.)  Then the idea is to let
$g$ be the so called increasing rearrangement of $f$, that is,
$$ g(x) = \sup\{t \in \R : \Pr(f < t) < x \} .$$
That $g$ has the required properties is easy to show.

Next, suppose that $f$ is nowhere constant, that is, we have
that $\Pr(f = s) = 0$ for all $s \in S$.  Then it may be seen that
$g$ is a strictly increasing function.  In that case, it may be seen
that the minimal complete sigma field for which $g$ is measurable
is the collection of Lebesgue measurable sets.

Now let us move on to state the main result.
We will start by setting our notation.
We will work on two probability spaces: a generic one 
$(\Omega,\cF,\Pr)$, and $([0,1]^\N,\cL_\N,\lambda)$.  
Here
$\cL_\N$ refers to the Lebesgue measurable sets on $[0,1]^\N$, and $\lambda$
refers to Lebesgue measure on $[0,1]^\N$.

Let us make some notational abuses.  The reason for this is to
make some expressions less cumbersome, while hopefully not
being too obscure.
We will always identify $[0,1]^n$ with the natural projection
of $[0,1]^\N$ onto the first $n$ coordinates.  
Any function $g$ on $[0,1]^n$ will be identified with its
canonical lifting on $[0,1]^\N$.  The notation $\cL_n$ refers
to the Lebesgue measurable sets on $[0,1]^n$, and also to their 
canonical lifting
onto $[0,1]^\N$. We will let $\lambda$ also refer
to Lebesgue measure on $[0,1]^n$.

Given a random variable $f$, and a sigma field $\cG$, we
will say that $f$ is {\it nowhere constant\/} with respect 
to $\cG$ if $\Pr(f = g) = 0$ for every $\cG$ measurable function
$g$.  Let us illustrate this notion with respect to a measurable function
$f$ on $[0,1]^2$.  It is nowhere constant with respect to the trivial
sigma field if and only if 
it is nowhere constant as defined above.  Let $\cG$ be the 
sigma field generated by the first coordinate.
Then $f$ is nowhere constant with respect to $\cG$ if and only if
for almost every $x \in [0,1]$ the function $y \mapsto f(x,y)$
is nowhere constant.  Now let $\cG$ be the set of all Lebesgue measurable
sets on $[0,1]^2$.  In this case, $f$ can never be nowhere constant with respect to
$\cG$.

If $f_1,\dots,f_n$ are random variables taking values in
a measurable space $(S,\cS)$, we will let 
$\sigma(f_1,\dots,f_n)$ denote the minimal sigma
field which contains all sets of measure zero,
and for which $f_1,\dots,f_n$ are measurable.
Note that a simple argument shows that
the $\sigma(f_1,\dots,f_n)$ measurable functions
coincide precisely with the functions that are
almost everywhere equal to
$F(f_1,\dots,f_n)$, where $F$ is some measurable
function on $S^n$.

Throughout these proofs we will use the following idea many times.
Let
$(\Omega,\cF,\Pr)$ be a 
probability space, and let $(S,\cS)$ be a
measurable space.  Suppose that $X:\Omega\to S$ is a measurable function
such that
$\sigma(X)$ is the whole
of $\cF$.  If $Y$ is any random variable on $\Omega$ taking values in $\R$, 
then
there exists a measurable map $\phi:S\to\R$ such that $Y = \phi \circ X$
almost surely.  This is easily seen by approximating $Y$ by simple functions.

\proclaim Theorem 2.1.  Let $(f_n)$ be a sequence of random variables
taking values in a separable measurable space $(S,\cS)$.
Then there exists a sequence of measurable functions
$(g_n:[0,1]^n \to S)$ that has the same law as $(f_n)$.
\moreproclaim
If further we have that $f_{n+1}$ is nowhere constant with respect to
$\sigma(f_1,\dots,f_n)$ for all $n \ge 0$,
then we may suppose that $\sigma(g_1,\dots,g_n) = \cL_n$, for all $n \ge 1$.
\moreproclaim
If even further we have that $(S,\cS)$ is $\R$ with the Borel sets, then
we may suppose $g_n$ is Borel measurable (that is, the preimages
of Borel sets are Borel sets),
and that $g_n(x_1,\dots,x_n)$ is a strictly increasing
function of $x_n \in (0,1)$ for almost all $x_1,\dots,x_{n-1}$.

\Proof:  
We may suppose without loss of generality that $S = \R$, and
$\cS$ is the Borel subsets of $\R$.  To see that we may do
this, see first that
we may suppose without loss of generality that $\cS$ separates
points in $S$, that is, if $s\ne t \in S$, then there exists
$A \in \cS$ such that $s \in A$ and $t \notin A$.
Let $\{\cC_n\}$ be a countable collection of sets in $\cS$
that generate $\cS$.  
Notice then that the sequence $\{\cC_n\}$ also separates points in $S$,
that is, if $s\ne t \in S$, then there exists a number $n$ such that
one and only one of $s$ or $t$ is in $\cC_n$.
Define a map $\varphi:S \to \R$:
$$ \varphi(s) = \sum_{m=1}^\infty {I_{(s \in C_m)}\over 3^m} .$$
Clearly $\varphi$ is injective.  Further
$\varphi$ maps $C_n$ to $D_n \cap \varphi(S)$, where
$$ D_n = \left\{ \sum_{m=1}^\infty {e_m \over 3^m} : e_m = 0 \hbox{ or }
   1, \, e_n = 0 \right\} , $$
and thus $\varphi$ maps any element of $\cS$ to a Borel subset of $\R$
intersected with $\varphi(S)$.  Conversely, the preimage of $D_n$ under
$\varphi$ is $C_n$, and hence the preimage of any Borel set under $\varphi$
is in $\cS$.  (This argument may be found in Chapter~1 of [DM]).

Now apply the theorem to $(\varphi \circ f_n)$, to obtain $(g_n)$.  Since
the law of $g_n$ is the same as $\varphi \circ f_n$, the range of
$g_n$ lies in the range of $\varphi(S)$ with probability one.  
Then the sequence $(\varphi^{-1} \circ g_n)$ will have the same
law as $(f_n)$.

So let us suppose as the induction hypothesis that we have
obtained $(g_n:[0,1]^n \to \R)_{1 \le n \le N}$ that has the same law as
$(f_n)_{1 \le n \le N}$, and that $g_n$ is Borel measurable,
for all $n \le N$, 
where $N$ is a non-negative integer.
(The induction is started with $N=0$ in which case the hypothesis
is vacuously true.  The arguments that follow simplify greatly in this case.)

For each $t \in \Q$, let 
$$ p_t = \E(I_{(f_{N+1}<t)} | \sigma(f_1,\dots,f_N)) .$$
Since $p_t$ is $\sigma(f_1,\dots,f_N)$ measurable, we may write
$$ p_t = q_t(f_1,\dots,f_N) \as$$
for some Borel measurable function $q_t:\R^N \to [0,1]$.  Define the 
Borel measurable functions $r_t$ by
$$ r_t = q_t(g_1,\dots,g_N) .$$
Since the sequence $(r_t)_{t \in \Q}$ has the same law as $(p_t)_{t \in \Q}$,
we see that there is a Borel set $B \subset [0,1]^N$ of full measure
such that if $(x_1,\dots,x_N) \in B$, then
$r_t (x_1,\dots,x_N)$ is an increasing function of $t \in \Q$, 
that tends to $0$ as $t \to -\infty$, and tends to $1$ as $t \to \infty$.

If $(x_1,\dots,x_N) \in B$ and $x_{N+1} \in (0,1)$, let
$$ g_{N+1}(x_1,\dots,x_{N+1}) =
   \sup \{ t\in\Q : r_t(x_1,\dots,x_N) < x_{N+1} \} ,$$
and let it be zero otherwise.
We see that $g_{N+1}$ is Borel measurable, since
$g_{N+1}$ on $B \times (0,1)$ is the supremum of countably many
Borel measurable functions $(s_t)_{t \in \Q}$ where
$$ s_t(x_1,\dots,x_{N+1}) = 
   \cases{ t &if $r_t(x_1,\dots,x_N) < x_{N+1} $\cr
           -\infty &otherwise.\cr } $$
It is easy to see that $g_{N+1}$ is 
always finite.

Let us now show that $(f_1,\dots,f_{N+1})$ has the same law as
$(g_1,\dots,g_{N+1})$.  It is sufficient to show that 
for any Borel set $A \subset \R^N$, and for
$t \in \Q$, that
$$ \Pr((f_1,\dots,f_N) \in A \hbox{ and } f_{N+1}<t)
   =
   \lambda((g_1,\dots,g_N) \in A \hbox{ and } g_{N+1}<t) .$$
So let us begin.
$$ \eqalignno{
   &\lambda((g_1,\dots,g_N) \in A \hbox{ and } g_{N+1}<t) \cr
   &=
   \lambda(\{(g_1,\dots,g_N) \in A \hbox{ and } g_{N+1}<t\} \cap
           B\times(0,1)) \cr
   &\le
   \lambda(\{(x_1,\dots,x_{N+1}):(g_1(x_1),\dots,g_N(x_1,\dots,x_N)) \in A 
           \hbox{ and } r_t(x_1,\dots,x_N) \ge x_{N+1} \} \cr
   &\hbox to 5.5 true in{}   
	   \cap B\times(0,1)) \cr
   &=
   \lambda(\{(x_1,\dots,x_{N+1}):(g_1(x_1),\dots,g_N(x_1,\dots,x_N)) \in A 
           \hbox{ and } r_t(x_1,\dots,x_N) \ge x_{N+1}\}) \cr
   &= \int_{[0,1]^N} \int_0^1 I_{(r_t \ge x)} I_{((g_1,\dots,g_N)\in A)}
      \, dx \, d\lambda \cr
   &=
   \int_{[0,1]^N} r_t I_{((g_1,\dots,g_N)\in A)} \, d\lambda \cr
   &=
   \E(p_t I_{((f_1,\dots,f_N)\in A)}) \cr
   &=
   \E(I_{(f_{N+1}<t)} I_{((f_1,\dots,f_N) \in A)}) \cr
   &=
   \Pr((f_1,\dots,f_N) \in A \hbox{ and } f_{N+1}<t) . \cr } $$
Similarly
$$ \lambda((g_1,\dots,g_N) \in A \hbox{ and } g_{N+1} \le t) 
   \ge
   \Pr((f_1,\dots,f_N) \in A \hbox{ and } f_{N+1}<t) , $$
and so by Lebesgue's monotone convergence theorem
$$ \eqalignno{
   \lambda((g_1,\dots,g_N) \in A \hbox{ and } g_{N+1} < t)
   &=
   \lim_{s \nearrow t\atop s\in\Q}
   \lambda((g_1,\dots,g_N) \in A \hbox{ and } g_{N+1} \le s) \cr
   &\ge
   \lim_{s \nearrow t\atop s\in\Q}
   \Pr((f_1,\dots,f_N) \in A \hbox{ and } f_{N+1} < s) \cr
   &=
   \Pr((f_1,\dots,f_N) \in A \hbox{ and } f_{N+1} < t) . \cr} $$
Thus we have shown that $(f_1,\dots,f_{N+1})$ has the same law as
$(g_1,\dots,g_{N+1})$.

Now let us add the assumption that 
$f_{n+1}$ is nowhere constant with respect to
$\sigma(f_1,\dots,f_n)$ for all $n \ge 0$.
Let us assume the inductive hypothesis that
$\sigma(g_1,\dots,g_N) = \cL_N$, 

Then there exist Borel measurable
functions
$(\alpha_n:\R^n \to [0,1])_{1\le n \le N}$ such that for $1 \le n \le N$
we have
$$ \alpha_n(g_1(x_1),\dots,g_N(x_1,\dots,x_N)) = x_n \ae.$$
Let $\beta_n = \alpha_n(g_1,\dots,g_N)$, and let $\gamma_n =
\alpha_n(f_1,\dots,f_N)$.

Notice that $g_{N+1}(x_1,\dots,x_{N+1})$ is an increasing function in 
$x_{N+1} \in (0,1)$
for fixed $(x_1,\dots,x_N) \in B$.  
Let us show that $g_{N+1}(x_1,\dots,x_{N+1})$ is a 
strictly increasing function in 
$x_{N+1}\in (0,1)$ for $(x_1,\dots,x_N) \in C$
where $C$ is a subset of $B$ of full measure.  For suppose otherwise.
Then the following set has positive measure:
the set of $(x_1,\dots,x_{N+1})$ such that for some rational number
$t \in [0,1]$ we have that $g_{N+1}(x_1,\dots,x_N,x_{N+1}) = 
g_{N+1}(x_1,\dots,x_N,t)$.  But this is equal to the set
$$ \bigcup_{t\in[0,1]\cup\Q} \{g_{N+1} 
   = g_{N+1}(\beta_1,\dots,\beta_N,t)\} .$$
But the measure of this set is the same as the probability of the set
$$ \bigcup_{t\in[0,1]\cup\Q} \{f_{N+1} 
   = g_{N+1}(\gamma_1,\dots,\gamma_N,t)\} .$$
Since $f_{N+1}$ is nowhere constant with respect to $\sigma(f_1,\dots,f_N)$,
it follows that this last set has probability zero.

Let us now show that $\sigma(g_1,\dots,g_{N+1}) = \cL_{N+1}$.
For any $a \in [0,1]$
we have that
$$ \{ g_{N+1} \le g_{N+1}(\beta_1,\dots,\beta_N,a) \} 
   = [0,1]^N \times [0,a] \quad\hbox{up to a set of measure zero}.$$
Thus $\sigma(g_1,\dots,g_{N+1})$ contains all set of the form
$[0,1]^N \times [0,a]$.  Since it also contains $\sigma(g_1,\dots,g_N)
= \cL_N$, the result follows.
\QED

\beginsection 3.\ \ Representation of martingales --- 1

\proclaim Theorem 3.1.  Let $(d_n)$ be a Bochner integrable 
martingale difference sequence
taking values in a Banach space $X$.
Then there exists a sequence of Bochner measurable functions
$(e_n:[0,1]^n \to X)$ such that $(d_n)$ has the same law as
$(e_n)$, and such that for almost every $x_1,\dots,x_n$ we have
$$ \int_0^1 e_n(x_1,\dots,x_n) \, dx_n = 0 .$$

\Proof:  
First we show that without loss of generality that
we may assume that $d_{n+1}$ is nowhere constant with respect
to $\sigma(d_1,\dots,d_n)$ for all $n \ge 0$.  To do this, replace
$\Omega$ by $\Omega\times[0,1]^\N$, replace $X$ by $X \times \R$
(where $\R$ is equipped with the sigma field of Borel sets),
and replace $d_n$ by $(\omega,(x_n)) \mapsto (d_n(\omega),x_n-{1\over2})$.
Apply the theorem to this sequence.  After obtaining the resulting
$(e_n)$, compose these functions with the natural projection of
$X\times\R$ to $X$.

Further, since $(d_n)$ is Bochner integrable, we may suppose without loss of
generality that $X$ is separable.

Apply Theorem~2.1 to obtain the sequence $(e_n)$.  Suppose that
$\phi:[0,1]^{n-1} \to \R$ is any bounded
measurable function.  Then
there exists a bounded Borel measurable
function $\theta:X^{n-1}\to\R$ such that
$\phi(x_1,\dots,x_{n-1}) = \theta(e_1(x_1),\dots,e_{n-1}(x_1,\dots,x_{n-1}))$
for almost all $x_1,\dots,x_{n-1}$.
Then
$$ \int_{[0,1]^n} \phi e_n \, d\lambda
   =
   \E(\theta(d_1,\dots,d_{n-1}) d_n) ,$$
and this is zero because $\E(d_n | \sigma(d_1,\dots,d_{n-1})) = 0$.
The result follows.
\QED

\beginsection 4.\ \ Representation of tangent sequences

The following definition may be found in [KW].
Let $(\cF_n)_{n\ge 0}$ be an 
increasing sequence of sigma fields on $(\Omega,\cF,
\Pr)$, where $\cF_0$ is the trivial sigma field.
Two adapted sequences $(f_n)$ and $(g_n)$ taking values in a measurable space
$(S,\cS)$
are said to be {\it
tangent\/} if for each $n \ge 1$ we have that the law of $f_n$ conditionally
on $\cF_{n-1}$ is the same as the law of $g_n$ conditionally on $\cF_{n-1}$,
that is, $\E(I_{(f_n \in A)}|\cF_{n-1}) = \E(I_{(g_n \in A)}|\cF_{n-1})$
for any $A \in \cS$.

Notice then that if $F:S^{2n-1} \to \R$ is any measurable bounded function,
that
$$ \E(F(f_1,g_1,\dots,f_{n-1},g_{n-1},f_n))
   =
   \E(F(f_1,g_1,\dots,f_{n-1},g_{n-1},g_n)) .$$
This is easily seen by reducing it to the case where $F$ is 
a finite linear combination of characteristic
functions of sets
of the form $A_1\times\dots\times A_{2n-1}$ where $A_1,\dots,A_{2n-1}
\in \cS$.

\proclaim Theorem 4.1.  Let $(S,\cS)$ be a separable measurable space.
Let $(f_n)$ and $(g_n)$ be 
$S$-valued measurable sequences adapted to
$(\cF_n)_{n \ge 0}$ that are tangent.  
Then there exists a sequence $(h_n:[0,1]^n \to 
S)$, and a sequence of Borel
measurable functions $(\phi_n:[0,1]^n \to [0,1])$ such
that the map $\phi_n(x_1,\dots,x_{n-1},\cdot)$ is a measure preserving
map on $[0,1]$ for almost every $x_1,\dots,x_{n-1}$,
such that the law of $(f_n,g_n)$ is the same as the law of
$(h_n,k_n)$, where 
$$ k_n(x_1,\dots,x_n) = h_n(x_1,\dots,x_{n-1},\phi_n(x_1,\dots,x_n)) .$$
Further, if $S$ is a separable Banach space, and
$(f_n)$ is a martingale difference sequence 
with respect to $\cF_n$, then
$$ \int_0^1 h_n(x_1,\dots,x_n) \, dx_n = 0 $$
for almost all $x_1,\dots,x_{n-1}$.

\proclaim Lemma 4.2.  Let $(S,\cS,\mu)$, $(T,\cT,\nu)$ 
and $(U,\cU,\tau)$ be  
measure spaces such that $\cU$ is separable.
Let $\phi:S\times T \to U$ be a measurable map such that for $A \in \cS$
and $C \in \cU$ we have that
$$ \mu \otimes \nu((s,t):s \in A \hbox{ and } \phi(s,t) \in C)
   = \mu(A) \tau(C) .$$
Then $\phi(s,\cdot)$ is a measure preserving map $T \to U$ for $\mu$
almost every $s$.

\Proof:
The hypothesis of the lemma can be recast as:
$$ \int_{s \in A} \nu(t:\phi(s,t) \in C) \, d\mu(s)
   =
   \int_A \tau(C) \, d\mu ,$$
from which we conclude that for every $C \in \cU$, we have that
$\nu(t:\phi(s,t) \in C) = \tau(C)$ for $\mu$ almost every $s$.

Let $\cC$ be a countable subcollection of $\cU$ that generates $\cU$.  Then
there exists a set $A \in \cS$ of full measure such that
$\nu(t:\phi(s,t) \in C) = \tau(C)$ for $C \in \cC$ and $s\in A$.
Let
$$ \cM = \{C \in \cU : 
\hbox{$\nu(t:\phi(s,t) \in C) = \tau(C)$ for $s\in A$}\} .$$
Notice that $\cM$ is a sigma field that contains $\cC$, and hence
$\cM = \cU$.
\QED

\Proofof Theorem~4.1:  
Without loss of generality (at least until we prove the last
part), let us suppose that $S = \R$ with the Borel sets.  To do this,
we use the same argument that we gave at the beginning of the proof of 
Theorem~2.1.

Consider the sequence $(f_1,g_1,f_2,g_2,\dots)$.
By the kind of argument we presented at the beginning
of Theorem~3.1, we may suppose
that each element of this sequence is nowhere constant with respect to
the sigma field generated by those elements in the sequence that precede
it.

Applying Theorem~2.1, we create a sequence of Borel measurable functions 
$(\tilde h_1, \tilde k_1, \tilde h_2, \tilde k_2,\dots)$ with the same law,
where $\tilde h_n:[0,1]^{2n-1} \to \R$, $\tilde k_n:[0,1]^{2n} \to \R$, and
$\tilde h_n(x_1,y_1,\dots,x_{n-1},y_{n-1},x_n)$ is a strictly increasing
function of $x_n\in(0,1)$ for almost every $x_1,y_1,\dots,x_{n-1},y_{n-1}$.

Define the function
$$ r_n(x_1,y_1,\dots,x_{n-1},y_{n-1},s)
   =
   \sup \{ t\in\Q \cap (0,1) : 
   \tilde h_n(x_1,y_1,\dots,x_{n-1},y_{n-1},t) < s \}  ,$$
where we will set the first supremum equal to zero if the set is empty.
Notice that $r_n$ is the supremum of countably many Borel measurable functions,
and hence is Borel measurable.  Also, since 
$\tilde h_n(x_1,y_1,\dots,y_{n-1},x_{n-1},x_n)$ is a strictly increasing
function of $x_n \in (0,1)$ 
for almost every $x_1,y_1,\dots,x_{n-1},y_{n-1}$,
we see that
$$ r_n(x_1,y_1,\dots,x_{n-1},y_{n-1},\tilde h_n(x_1,y_1,\dots,x_{n-1},
       y_{n-1},x_n)) = x_n $$
for almost every $x_1,y_1,\dots,x_n$.
Define
$$ \tilde \phi_n(x_1,y_1,\dots,x_n,y_n) = 
   r_n(x_1,y_1,\dots,x_{n-1},y_{n-1},
   \tilde k_n(x_1,y_1,\dots,x_n,y_n)) .$$
Now let $\alpha_n,\beta_n$ be functions such that
$$ \eqalignno{
   x_n &= \alpha_n
          (\tilde h_1(x_1),\tilde k_1(x_1,y_1),\dots,
	  \tilde h_n(x_1,y_1,\dots,x_n)) \ae \cr
   y_n &= \beta_n
          (\tilde h_1(x_1),\tilde k_1(x_1,y_1),\dots,
	  \tilde k_n(x_1,y_1,\dots,x_n,y_n))
   \ae ,\cr}$$
and let $\gamma_n,\delta_n$ be random variables defined by
$$ \eqalignno{
   \gamma_n &= \alpha_n(f_1,g_1,\dots,f_n) \cr
   \delta_n &= \beta_n(f_1,g_1,\dots,f_n,g_n) . \cr } $$

If $A$ is any Borel subset of $[0,1]^{2n-2}$, and $C$ is a Borel
subset of $[0,1]$, then
since $f_n$ and $g_n$ have the same law with respect to $\cF_{n-1}$,
$$ \eqalignno{
   &\lambda((x_1,y_1,\dots,x_{n-1},y_{n-1})\in A \hbox{ and }
            \tilde \phi_n(x_1,y_1,\dots,x_n,y_n) \in C) \cr 
   &=
   \Pr((\gamma_1,\delta_1,\dots,\gamma_{n-1},\delta_{n-1}) \in A \hbox{ and }
       r_n(\gamma_1,\delta_1,\dots,\gamma_{n-1},\delta_{n-1},
       g_n) \in C) \cr 
   &=
   \Pr((\gamma_1,\delta_1,\dots,\gamma_{n-1},\delta_{n-1}) \in A \hbox{ and }
       r_n(\gamma_1,\delta_1,\dots,\gamma_{n-1},\delta_{n-1},
       f_n) \in C) \cr
   &=
   \lambda((x_1,y_1,\dots,x_{n-1},y_{n-1}) \in A \hbox{ and }
       r_n(x_1,y_1,\dots,x_{n-1},y_{n-1},
       \tilde h_n(x_1,y_1,\dots,x_n)) \in C) \cr
   &=
   \lambda((x_1,y_1,\dots,x_{n-1},y_{n-1}) \in A \hbox{ and } x_n \in C) \cr
   &=
   \lambda(A) \lambda(C) . \cr } $$
Thus, by Lemma~4.2, we have
that
$\tilde \phi(x_1,y_1,\dots,x_{n-1},y_{n-1},\cdot,\cdot)$ is a measure
preserving map from $[0,1]^2$ to $[0,1]$
for almost every $x_1,y_1,\dots,x_{n-1},y_{n-1}$.

Now let $x \mapsto (\theta_1(x),\theta_2(x))$ be any Borel measurable
measure preserving
map from $[0,1]$ to $[0,1]^2$.  Let
$$ \eqalignno{
   h_n(x_1,\dots,x_n) 
   &= 
   \tilde h_n(\theta_1(x_1),\theta_2(x_1), \dots, \theta_1(x_n)) \cr
   k_n(x_1,\dots,x_n)
   &=
   \tilde k_n(\theta_1(x_1),\theta_2(x_1),\dots,\theta_1(x_n),\theta_2(x_n))
   .\cr} $$
Notice that
$$ k_n(x_1,\dots,x_n) = h_n(x_1,\dots,x_{n-1},
   \phi_n(x_1,\dots,x_n)) ,$$
where
$$ \phi_n(x_1,\dots,x_n) =
   \tilde\phi_n(\theta_1(x_1),\theta_2(x_1),\dots,\theta_1(x_n),\theta_2(x_n))
   .$$
The last part of the theorem follows by exactly the same argument as
for Theorem~3.1.
\QED

\beginsection 5.\ \ Representation of sequences satisfying condition (C.I.)

The following definition may be found in [KW].
Let $(\cF_n)_{n\ge 0}$ be an increasing sequence of sigma fields on 
$(\Omega,\cF,
\Pr)$, where $\cF_0$ is the trivial sigma field.
An adapted sequence $(f_n)$ taking values in a measurable space
$(S,\cS)$ is said to satisfy {\it condition (C.I.)\/}
if there exists a sigma field $\cG \subset \cF$ such that 
the law of $f_n$ conditionally on $\cF_{n-1}$ is the same
as the law of $f_n$ conditionally on $\cG$, that is,
$\E(I_{(f_n \in A)}|\cF_{n-1}) = \E(I_{(f_n \in A)}|\cG)$
for any $A \in \cS$, and if the sequence $(f_n)$ is conditionally
independent with respect to $\cG$, that is, for any sequence of sets
$A_n \in \cS$ we have
$$ \E(I_{(f_1 \in A_1)} \cdots I_{(f_n \in A_n)} | \cG)
   =
   \E(I_{(f_1 \in A_1)} | \cG) \cdots \E(I_{(f_n \in A_n)} | \cG) .$$

It is shown in [KW] that given any sequence $(f_n)$ adapted to some
filtration, that after possibly enlarging the underlying probability
space, that there exists a sequence $(\tilde f_n)$ that is tangent to
$(f_n)$ and that satisfies condition~(C.I.).  

Indeed, using Theorem~2.1,
we can show a technically weaker but essentially identical result as follows.
Let $(g_n)$ be the sequence constructed by Theorem~2.1.  Enlarge the
probability space $[0,1]^\N$ to $[0,1]^\N \times [0,1]^\N$, and define
the sequences
$$ \eqalignno{
   u_n((x_n),(y_n)) &= g_n(x_1,\dots,x_{n-1},x_n) \cr
   v_n((x_n),(y_n)) &= g_n(x_1,\dots,x_{n-1},y_n) .\cr } $$
Then $(u_n)$ has the same law as $(f_n)$, and $(v_n)$ is tangent to
$(u_n)$ with respect to the filtration $(\cL_n \otimes \cL_0)$, and
$(v_n)$ satisfies condition~(C.I.).
(Here $\cL_n$ is the minimal complete sigma field on $[0,1]^\N$
for which the first $n$ coordinate functions are measurable.  In particular,
$\cL_0$ denotes the trivial complete measure sigma field on $[0,1]^\N$.)

This is the motivation behind the next result, which shows that this
previous construction is the canonical representation.

\proclaim Theorem 5.1.  Let $(S,\cS)$ be a separable measurable space.
Let $(f_n)$ be an $S$-valued measurable sequence
adapted to
$(\cF_n)_{n \ge 0}$ that satisfies condition (C.I.).
Then there exists a sequence $(h_n:[0,1]^n \to 
S)$ such that if we define the functions
$u_n$ on $[0,1]^\N \times [0,1]^\N$ by
$$ u_n((x_n),(y_n)) = h_n(x_1,\dots,x_{n-1},y_n) ,$$
then the laws of $(f_n)$, $(h_n)$ and $(u_n)$ are the same.
\moreproclaim
Further, if $S$ is a separable Banach space, and
$(f_n)$ is a martingale difference sequence 
with respect to $\cF_n$, then
$$ \int_0^1 h_n(x_1,\dots,x_n) \, dx_n = 0 $$
for almost all $x_1,\dots,x_n$.

\Proof:  As usual, we suppose without loss of generality that
$f_n$ is nowhere constant with respect to $\sigma(f_1,\dots,f_{n-1})$.
Let $(h_n)$ be the sequence generated by Theorem~2.1.
Then for any $A_1,\dots,A_n \in \cS$ we have that
$$ \eqalignno{
   &\Pr((f_1,\dots,f_n) \in A_1 \times \cdots \times A_n) \cr
   &=
   \E(\E(I_{(f_1 \in A_1)} \cdots I_{(f_n \in A_n)} | \cG)) \cr
   &=
   \E(\E(I_{(f_1 \in A_1)} | \cG) \cdots \E(I_{(f_n \in A_n)} | \cG)) \cr
   &=
   \E(\E(I_{(f_1 \in A_1)} | \cF_0) \cdots 
   \E(I_{(f_n \in A_n)} | \cF_{n-1})) .\cr}$$
Now for each $k=1,\dots,n$,
let $F_k:S^{k-1} \to [0,1]$ be the measurable function such that
$$ \E(I_{(f_k \in A_k)} | \cF_{k-1}) = F_k(f_1,\dots,f_{k-1}) \as.$$
Notice that if $G:S^{n-1} \to \R$ is bounded measurable, then
$$ \E(I_{(f_k \in A_k)} G(f_1,\dots,f_{k-1}))
   =
   \E(F_k(f_1,\dots,f_{k-1}) G(f_1,\dots,f_{k-1})) ,$$
and hence
$$ \E(I_{(h_k \in A_k)} G(h_1,\dots,h_{k-1}))
   =
   \E(F_k(h_1,\dots,h_{k-1}) G(h_1,\dots,h_{k-1})) .$$
Since $\sigma(h_1,\dots,h_{k-1}) = \cL_{k-1}$, we see that
$$ \E(I_{(h_k \in A_k)} | \cL_{k-1}) = F_k(h_1,\dots,h_{k-1}) \ae. $$
Thus the sequence 
$$ (\E(I_{(f_k \in A_k)} | \cF_{k-1}))_{1 \le k \le n} $$
has the same law as the sequence
$$ (\E(I_{(h_k \in A_k)} | \cL_{k-1}))_{1 \le k \le n} ,$$
and this last sequence is very easily seen to have
the same law as the sequence
$$ (\E(I_{(u_k \in A_k)} | \cL_{k-1}\otimes\cL_0))_{1 \le k \le n} .$$
Hence
$$ \eqalignno{
   &\E(\E(I_{(f_1 \in A_1)} | \cF_0) \cdots 
   \E(I_{(f_n \in A_n)} | \cF_{n-1})) \cr
   &=
   \int(\E(I_{(u_1 \in A_1)} | \cL_0\otimes\cL_0) \cdots 
   \E I_{(u_n \in A_n)} | \cL_{n-1}\otimes\cL_0) \, 
   d \lambda\otimes\lambda ,\cr}$$
and this last quantity is easily computed to be
$$ \lambda\otimes\lambda((u_1,\dots,u_n) \in A_1 \times \cdots \times A_n) .$$
Finally, the last part of the theorem follows exactly as in Theorem~3.1.
\QED

\beginsection 6.\ \ Tangent sequences in UMD spaces

The concept of UMD spaces was introduced by Aldous [A],
and extensively explored by Burkholder [Bu].  It is from this
second reference that
the following definition of UMD (along with many other
equivalent definitions) may be found.

A Banach space $X$ is said to be {\it UMD\/} if for some
$1<p<\infty$ (equivalently all $1<p<\infty$), there exists
a positive constant $c$ depending only upon $X$ and $p$ such that
if $(d_n)$ is a Bochner integrable $X$-valued martingale difference
sequence, then for any signs $\epsilon_n = \pm 1$, we have that
$$ \normo{ \sum_{n=1}^N \epsilon_n d_n}_{L_p(X)}
   \le
   c
   \normo{ \sum_{n=1}^N d_n }_{L_p(X)} .$$
This section is devoted to providing a new proof of the following result.
This is essentially the same as the first part of Theorem~2.2 from the
paper by McConnell [M2].  
While the following proof
can be rewritten so as to avoid needing the representation theorems,
to do so would be to remove the motivation behind this proof.

\proclaim Theorem 6.1.  Let $X$ be a UMD space.  Then given $1<p<\infty$,
there is a positive constant $C$, depending only upon $X$ and $p$, such that
given a filtration $(\cF_n)_{n \ge 0}$, 
and two martingales $(f_n)$ and $(g_n)$ adapted
to this filtration
that are $X$-valued Bochner integrable and tangent to each other,
we have that
$$ \normo{\sum_{n=1}^N g_n}_{L_p(X)}
   \le C
   \normo{\sum_{n=1}^N f_n}_{L_p(X)} .$$

\Proof:  It is sufficient to prove the inequality in the case when one
of the sequences satisfies condition~(C.I.), since the inequality
may be then obtained by comparing both sides with another sequence
tangent to them both that satisfies condition~(C.I.).  We will
prove it in the case when the sequence $(g_n)$ satisfies condition~(C.I.),
but it will be evident that the proof also works when the sequence
$(f_n)$ satisfies condition~(C.I.).

So applying Theorems~4.1 and~5.1, we see that there exists a
sequence $h_n:[0,1]^n \to X$ of $X$-valued Bochner integrable functions
such that 
$$ \int_0^1 h_n(x_1,\dots,x_n) \, dx_n = 0 $$
for almost all $x_1,\dots,x_{n-1}$, and such that if we define the
the functions $d_n,e_n:[0,1]^\N \times [0,1]^\N \to X$ by
$$ \eqalignno{
   d_n((x_n),(y_n)) &= h_n(x_1,\dots,x_{n-1},x_n) \cr
   e_n((x_n),(y_n)) &= h_n(x_1,\dots,x_{n-1},y_n) ,\cr } $$
then the sequence $(f_1,g_1,f_2,g_2,\dots)$ has the same law as
the sequence $(d_1,e_1,d_2,e_2, \dots)$.

Now we define a martingale difference sequence 
$r_n$ on $[0,1]^\N \times [0,1]^\N$ by
$ r_{2n-1} = d_n + e_n $
and
$ r_{2n} = d_n - e_n $,
where the filtration is $\cG_{2n} = \cL_n \otimes \cL_n$, and
$\cG_{2n-1}$ is the collection of those sets in $\cG_{2n}$ generated
by sets like 
$$ \{((x_n),(y_n)) : (x_1,\dots,x_n) \in A\times C,\, 
                     (y_1,\dots,y_n) \in B \times C \} ,$$
where $A,B \in \cL_{n-1}$,
and $C \in \cL_1$.

Then we see that 
$ \sum_{n=1}^N d_n = {1\over 2} \sum_{n=1}^{2N} r_n $
and
$ \sum_{n=1}^N e_n = {1\over 2} \sum_{n=1}^{2N} (-1)^{n+1} r_n  $,
and thus the result follows immediately from the fact that $X$ is a UMD space.
\QED

\beginsection 7.\ \ Representations of martingales --- 2

In this section, we present a version of a result of McConnell [M1].
This result provides a Skorohod like representation for vector
valued martingales.  Our result is weaker than that of McConnell in that
his process is adapted to a filtration generated by a one-dimensional
Brownian motion.  The process presented here is adapted to a two
dimensional Brownian motion.

\proclaim Theorem 7.1.  Let $(d_n)$ be a Bochner integrable
martingale difference sequence
taking values in a Banach space $X$.
Then there exists a continuous time $X$-valued stochastic process
$(F_t)_{t \ge 0}$ with continuous sample paths
such that $(d_n)$ has the same law as
$(F_n - F_{n-1})$.  Furthermore, the process $(F_t)$ is adapted to
a two dimensional brownian motion.

\Proof:  Let $D = \{z\in \C : \modo z < 1\}$, 
$\bar D = \{z\in \C : \modo z \le 1\}$
and $\partial D = \{z\in \C : \modo z = 1\}$, and give $\partial D$ normalized
Lebesgue measure $d\modo{z} / 2\pi$.
Since $\partial D$ is measure equivalent to $[0,1]$, we
may apply Theorem~3.1 to obtain a sequence $(e_n:\partial D^n \to X)$ that
has the same law as $(f_n)$, and such that 
$$ \int_{\partial D} e_n(z_1,\dots,z_n) \, {d\modo{z_n} \over 2 \pi} = 0 $$
for almost every $z_1,\dots,z_{n-1}$.

For each $z_1,\dots,z_{n-1}$, extend the function $z_n \mapsto
e_n(z_1,\dots,z_n)$ to its harmonic extension on $\bar D$.  We
will identify $e_n$ with this extension.

Now let $(b_t^{(n)})_{n\in\N}$ be a sequence of independent Brownian
motions into $\C$ with origin at zero.  Let $\tau_n$ be the stopping
time $\tau_n = \inf\{t:\smodo{b_t^{(n)}} \ge 1 \}$.  Notice that
$\tau_n$ is finite almost surely.

By the It\^o calculus, we have that
$$ e_n(b^{(1)}_{\tau_1},\dots,b^{(n)}_{\min\{t,\tau_n\}}) =
   \int_0^{\min\{t,\tau_n\}} \nabla_n e_n
   (b^{(1)}_{\tau_1},\dots,b^{(n-1)}_{\tau_{n-1}},b^{(n)}_s) \cdot 
   db_s^{(n)} ,$$
where $\nabla_n$ refers to taking the gradient with respect to the
$n$th coordinate.
Define a process
$$ G_{n,t} = \left(\sum_{m=1}^{n-1}  
   e_m(b^{(1)}_{\tau_1},\dots,b^{(m)}_{\tau_m}) \right) +
   e_n(b^{(1)}_{\tau_1},\dots,b^{(n)}_{\min\{t,\tau_n\}}) .$$
Let $\phi:[0,1) \to [0,\infty)$ be defined by
$\phi(t) = t/(1-t)$, and let
$$ F_t = G_{[t+1],\phi(t-[t])} ,$$
where $[t]$ represents the largest integer less than $t$.

To see that $(F_t)$ may be adapted to a two-dimensional brownian motion,
let $\tilde b_t$ be a brownian motion taking values in $\C$, and set
$$ b_t^{(n)} = \int_{n-1}^{n-1+\phi^{-1}(t)} \sqrt{\phi'(s)} \, d\tilde b_s .$$
\QED

\beginsection Acknowledgments

I would like to mention the help of Gerald Edgar and P. Fitzsimmons in
obtaining the argument at the beginning of Theorem~2.1, where it is
shown that without loss of generality that $(S,\cS)$ is $\R$ with
the Borel sets, and for bringing the reference [DM] to my attention.
I also extend thanks to Donald Burkholder, who brought the reference
[M1] to may attention.  Finally, I would like to thank Terry McConnell
and an anonymous referee
for many useful remarks.

\beginsection References

\frenchspacing

\item{[A]} D.J.~Aldous, Unconditional bases and martingales in 
$L\sb{p}(F)$, {\sl Math. Proc. Cambridge Philos. Soc. {\bf 85}, (1979),
117--123.}

\item{[Bu]} D.L.~Burkholder, A geometrical characterization of Banach 
spaces in which martingale difference sequences are unconditional,
{\sl Ann. Probab. {\bf 9}, (1981), 997--1011.} 

\item{[DM]} C.~Dellacherie and P.-A. Meyer, {\sl Probabilities and
Potential}, North-Holland Mathematics Studies 29, North-Holland,
Amsterdam-New York-Oxford, 1978.

\item{[KW]} S.~Kwapie\'n and W.A.~Woyczy\'nski, Tangent sequences of
random variables, in {\sl Almost Everywhere Convergence}, G.A.~Edgar
and L.~Sucheston, Eds., Academic Press, 1989, pp. 237--265. 

\item{[M1]} T.R.~McConnell, A Skorohod-like representation in infinite
dimensions, {\sl Probability in Banach spaces, V (Medford, Mass.,
1984), 359--368}, Lecture Notes in Math., 1153, Springer, Berlin-New York, 1985.

\item{[M2]} T.R.~McConnell, Decoupling and stochastic integration in 
UMD Banach spaces, {\sl Probab. Math. Statist. {\bf 10}, (1989),
283--295.}

\bye